\documentclass[12pt,reqno]{amsart}
\usepackage{amsmath,amssymb,eucal}
\usepackage[initials]{amsrefs}
\usepackage{bbm}
\usepackage{amscd}
\font\cmssl=cmss10 at 12 pt

\newcommand{\G}{\Gamma}

\newcommand{\bC}{\mathbb{C}}
\newcommand{\bR}{\mathbb{R}}

\newcommand{\id} {\mathbbm{1}}

\renewcommand{\square}{\kern1pt\vbox
               {\hrule height 0.6pt\hbox{\vrule width 0.6pt\hskip 3pt
    \vbox{\vskip 6pt}\hskip 3pt\vrule width 0.6pt}\hrule height0.6pt}
    \kern1pt}

\newcommand{\ra}{\rightarrow}

\newtheorem{Th}{Theorem}
\newtheorem{Prop}{Proposition}
\newtheorem{Cor}{Corollary}
\newtheorem{Lem}{Lemma}
\newtheorem{Def}{Definition}
\newtheorem{Rem}{Remark}
\newtheorem{Ex}{Example}
\newcommand{\bt}{\begin{Th}\ \ }
\newcommand{\et}{\end{Th}}
\newcommand{\bp}{\begin{Prop}\ \ }
\newcommand{\ep}{\end{Prop}}
\newcommand{\bc}{\begin{Cor}\ \ }
\newcommand{\ec}{\end{Cor}}
\newcommand{\bl}{\begin{Lem}\ \ }
\newcommand{\el}{\end{Lem}}
\newcommand{\bd}{\begin{Def}\ \ }
\newcommand{\ed}{\end{Def}}
\newcommand{\pf}{\begin{proof}}  
\newcommand{\epf}{\end{proof}}

\newcommand{\n}{\nabla}

\newcommand{\ot}{\otimes}

\newcommand{\be}{\begin{equation}}
\newcommand{\ee}{\end{equation}}

\newcommand{\arr}{\begin{array}{rlll}}
\newcommand{\ea}{\end{array}}
\newcommand{\bea}{\begin{eqnarray}}
\newcommand{\eea}{\end{eqnarray}}
\newcommand{\bean}{\begin{eqnarray*}}
\newcommand{\eean}{\end{eqnarray*}}

\catcode`@=11
\@addtoreset{equation}{section}
\catcode`@=12

\begin{document}
\title{Integrability of generalized pluriharmonic maps}
\author{Lars Sch\"afer}
\address{Leibniz Universit\"at Hannover, Institut f\"ur Differentialgeometrie,    
  Welfengarten 1, 30167 Hannover, Germany}
\email{schaefer@math.uni-hannover.de}
\date{\today}
\begin{abstract}
In this paper we provide examples of maps from almost complex domains into pseudo-Riemannian symmetric targets, which are pluriharmonic and not integrable, i.e. do not admit an associated family. More precisely, for one class of examples the source has a non-integrable complex structure, like for instance a nearly K\"ahler structure and the target is a Riemannian symmetric space and for the other class the source is a complex manifold and the target is a pseudo-Riemannian symmetric space. These examples show, that a former result, Theorem 5.3 of \cite{CS1}, on the existence of associated families is sharp.
\end{abstract}
\maketitle

\section{Introduction}
\noindent
It is known for a long time, that a minimal surface in 3-space (up to coverings) allows a one-parameter family of isometric deformations preserving principal curvatures and rotating principal curvature directions, called {\it the associated family.} The most famous example is the deformation of the catenoid to the helicoid. Existence of associated families has been established for harmonic maps from surfaces into compact symmetric spaces (cf. \cites{U,Hi,BFPP,DPW}) and affine symmetric spaces (cf. \cites{BD05,H98}). If one passes from Riemannian surfaces to K\"ahler manifolds, one has to consider {\it pluriharmonic maps}, i.e. maps, such that the restriction to arbitrary complex curves is harmonic. These pluriharmonic maps admit associated families, as was shown by \cite{OV} for maps into compact Lie groups and for Riemannian symmetric spaces $S=G/K$ in \cites{BFPP,ET}. In fact, pluriharmonic maps are in these cases characterized by the property of having an associated family. \\ 
In applications from theoretical physics, namely tt*-geometry \cites{S,S1,CS1}, it is necessary to consider non-integrable complex structures on the source and pseudo-Riemannian metrics on the target manifolds. For this general situation there are different ways of defining the notion of a pluriharmonic map: One  is to consider the partial differential equation $(\n d \phi)^{1,1}=0$ (see Definition \ref{pl_harm_equ_def} for details) and another to  demand the existence of an associated family (compare Definition \ref{S1_pl_def} for a precise statement). For maps 
$h \colon (M,J) \ra (P=G/K, J^P)$ into K\"ahler symmetric targets $(P, J^P)$ the authors of \cite{EQ} propose in Section 5 to define pluriharmonic maps as maps, such that the $\mathfrak g$-valued one-form $J^P \circ dh \circ J$ is closed, where $P$ is embedded by using the {\it standard embedding.} The first two notions are only equivalent under two additional assumptions (A1) and (A2) discussed in Theorem \ref{int_family}.
In the present article we construct pluriharmonic maps in the first sense without associated family each of these only violating either (A1) or (A2). At the same time this shows, that Theorem \ref{int_family} is sharp, in the sense, that one cannot omit the additional assumptions. \\
Let us shortly explain the structure of this paper: In the second section we recall definitions and results on harmonic maps, pluriharmonic maps and associated families. Further we introduce nearly K\"ahler manifolds and properties and examples of pluriharmonic maps from nearly K\"ahler sources admitting associated families. In the third section we give a general construction for pluriharmonic maps without associated families from non-integrable almost complex manifolds $(M,J).$ In particular, we obtain examples of such maps from $(M,J)$ into $S^{2n}$ and $\bC P^n$ and from nearly K\"ahler sources. In the last section we discuss a construction which yields pluriharmonic maps without associated families from complex manifolds into pseudo-Riemannian symmetric spaces. \\
\noindent
{\bf Acknowledgements.}
The author thanks F.\ E. Burstall for valuable discussions during a visit at the University of Bath.
\section{Harmonic and pluriharmonic maps}
\noindent 
In this section we recall some notions and results on {\it harmonic maps,  
pluriharmonic maps} and propose natural generalizations of the notion of a
pluriharmonic map.
\subsection{Harmonic maps}
Given two pseudo-Riemannian manifolds $(M,g)$ and $(N,h)$ one associates
to a smooth map $$ \phi \colon M \ra N$$ 
an {\cmssl energy density} $e(\phi)$ by $$e(\phi):= \langle d\phi,d\phi \rangle,$$ where we consider $d\phi$ as a section of $T^*M \ot \phi^*TN$ and
where $\langle \cdot ,\cdot \rangle$ is the pseudo-Riemannian bundle metric on $T^*M\ot \phi^*TN$ induced by $g$ and
$h.$ A map $ \phi \colon M \ra N$ is called {\cmssl harmonic} if it is a stationary point of the {\cmssl energy}
$$ E(\phi):= \int_D e(\phi) dv_g$$
of $\phi$ for all compact sub-domains $ D \subset M.$ Let us note, that in physics
harmonic maps are also called {\cmssl non-linear sigma
  models} \cites{BP,Gl}.   \\
The corresponding Euler-Lagrange equation or the {\cmssl harmonic map equation}
is \be \label{harm_equ} \mbox{Tr}_g \nabla d\phi =0, \ee
 where $\nabla$ is the connection on $T^*M\ot \phi^*TN$ induced by the
 Levi-Civita connections $\n^g$ and $\n^h$ of the metrics $g$ and $h.$
A natural generalization of the concept of a harmonic map \cites{H98,BD05} is to consider a map 
$$ \phi \colon M \ra N$$ into an affine manifold $(N,\n^N)$ and to replace $\nabla$ in equation \eqref{harm_equ} by the connection induced by $\n^N$ and $\n^g.$
\subsection{Pluriharmonic maps}
Let $(M,J)$ be a complex manifold and $(N,h)$ a pseudo-Riemannian
manifold.  The harmonic functional for a smooth map from a Riemannian surface $(\Sigma,g)$ into $(N,h)$ is conformally invariant. 
In consequence, the subsequent notion of a  pluriharmonic map does not depend on the metric chosen in the conformal class induced by the complex structure $J.$ \\
A map $\phi \colon M \ra N$ is called {\cmssl pluriharmonic}, if
its restriction $\phi \circ \iota \colon \Sigma \ra N$ to an arbitrary Riemannian surface $\iota \colon \Sigma \ra M$ is
harmonic. \\ 
Pluriharmonic maps are equivalently described by a partial differential equation: A map $\phi: M \ra N$ is  pluriharmonic if and only if it satisfies the equation 
\be (\nabla d\phi)^{1,1}=0, \label{pl_harm_equ} \ee
where $\nabla$ is the connection on $T^*M \ot \phi^*TN$ induced by a
torsion-free complex  connection $D,$ i.e. a torsion-free connection
satisfying $DJ=0,$ and by the Levi-Civita connection $\n^h$ of $h.$\\ 
\noindent 
One possible way to generalize pluriharmonic maps to maps from almost complex source manifolds is by considering maps satisfying the Equation \eqref{pl_harm_equ}. In order to do this, we have to fix a connection on $M.$ The existence of a torsion-free complex connection on $(M,J)$ is equivalent to the integrability of $J.$ As a consequence we need to weaken the requirements on the connection. This motivates to consider complex connections of {\it Nijenhuis type} in the sense of the next Definition. We shall see below, that the notion of a {\it generalized pluriharmonic map} does not depend on the choice of the Nijenhuis type complex connection.

\bd
Let $(M,J)$ be an almost complex manifold. A complex connection $D$ is called of {\cmssl
  Nijenhuis type} provided that its torsion $T$ satisfies the condition $T^{1,1}=0.$ 
\ed
\noindent 
The {\cmssl Nijenhuis tensor\footnote{In \cite{KN} the Nijenhuis tensor is
    defined with a factor 2.}} of the almost complex structure $J$ is defined as
  $$N_J(X,Y): =[JX,JY] -[X,Y] -J[X,JY] -J[JX,Y],\, X,Y \in \Gamma(TM).$$ 
A well-known result (see for example Theorem 3.4 in  ch. IX of \cite{KN}) ensures that on an almost complex manifold,
there exists a complex connection $D,$ such that its torsion $T$ is a multiple of the Nijenhuis-tensor.  As $N_J$ has  vanishing $(1,1)$-part,  the space of Nijenhuis type complex connections is not empty and is described as follows.
\bp \label{NC_conn_prop}
The space of Nijenhuis type complex connections is an affine space over 
$$\mathcal{NC}=\{S \in \Gamma\left( T^*M\ot \mbox{End}(TM,J) \right) \; | \;   \mbox{alt}(S)^{1,1}=0 \},$$ 
where 
$\mbox{End}(TM,J)$ are the $J$-linear endomorphisms and where 
$$\mbox{alt} \colon T^*M\ot \mbox{End}(TM) \ra T^*M\ot \mbox{End}(TM)$$ is the
alternation map  $$\mbox{alt}(S)(X,Y)= S_XY-S_YX, \, X,Y \in TM.$$
Moreover, one has 
\be \label{nice_con_two} \mathcal{NC} = \{S \in \Gamma\left( T^*M\ot \mbox{End}(TM,J) \right) \; | \;   S_Z \bar W = 0 \mbox{ for all } Z, W \in T^{1,0}M \}. \ee  
\ep
\pf
Let us consider two Nijenhuis type complex connections $D^1$ and $D^2$ and set 
$S_XY:= D^1_XY-D^2_XY$  for some vector fields $X,Y.$ It follows from $D^1 J=D^2 J=0,$ that one has 
$$0= (D^1_XJ) Y -(D^2_XJ) Y = [S_X,J]Y, \; \forall X,Y \in \Gamma(TM).$$ In other words $S_X$ is a complex linear endomorphism. The torsions $T^1$ and $T^2$ of $D^1$ and $D^2$ are related by
$$ T^1(X,Y)-T^2(X,Y)= S_XY-S_YX = \mbox{alt}(S)(X,Y).$$ 
Hence the condition, that the $(1,1)$-part of the torsion vanishes yields
$\mbox{alt}(S)^{1,1}=0,$  i.e. $ S_Z \bar W = S_{\bar W}Z$ for all $Z,W \in T^{1,0}M.$ \\ Conversely, let $D^1$ be a Nijenhuis type complex connection 
and $S \in \mathcal{NC},$ then $D^2=D^1+S$ is a complex connection, as $S_X$ is complex linear. Since the torsion of $D^2$ is $T^2=\mbox{alt}(S),$ $D^2$ is a connection of Nijenhuis type.\\
Let us now fix $S \in \mathcal{NC}$ and  $Z,W$ of type $(1,0).$ As $S_Z$ and $S_{\bar W}$ are complex linear, it follows, that $S_ Z \bar W$ is of type $(0,1)$ and $S_{\bar W} Z$ is of type $(1,0).$ From $S_ Z \bar W = S_{\bar W} Z$ we conclude $S_ Z \bar W=S_{\bar W} Z=0,$
which proves the inclusion '$\subset$' in equation  \eqref{nice_con_two}.
The other inclusion '$\supset$' in equation  \eqref{nice_con_two} is obvious. \\
\epf 

\bp \label{nice_conn_prop}
Let  $\phi \colon (M,J,D) \ra (N,h)$ be a map between an almost complex manifold  $(M,J)$  endowed with a complex connection $D$ of Nijenhuis type and a pseudo-Riemannian manifold $(N,h),$  then
\begin{itemize}
\item[(i)]
the field $(\nabla d\phi)^{1,1}$ is symmetric, i.e.  
$\nabla d\phi(Z,\bar W)=\nabla d\phi(\bar W,Z)$ for $Z,W \in T^{1,0}M$ 
\item[(ii)]
and the equation  $ (\nabla d\phi)^{1,1} =0$  does not depend on the choice of the Nijenhuis type complex connection $D.$
\end{itemize}
\ep
\pf In order to prove (i) we recall the following information: Given a smooth map $\phi \colon M \ra N,$  the differential  $\Phi:=d\phi :\, TM \ra \phi^*TN=E$ induces a vector bundle homomorphism between the tangent bundle of $M$ and the pull-back of $TN$ via $\phi.$ From the vanishing of the torsion of $\n^h$ we get the identity 
\bea
\n^E_V \Phi(W) - \n^E_W \Phi(V) - \Phi([V,W]) = 0, \mbox{ for } V,W \in \Gamma(TM),\label{Int_equ_I}
\eea
where $\n^E= \phi^*\n^h$ denotes the pull-back connection, i.e. the connection which is induced on $E$ by the Levi-Civita connection $\n^h$ of $h.$ 
 \\
For $Z,W$ of type $(1,0)$ we get using $T^{1,1}=0$ and equation \eqref{Int_equ_I}, that it holds
\bean
0&=& \n^E_Z \Phi(\bar W) - \n^E_{\bar W}\Phi(Z) - \Phi(D_Z \bar W - D_{\bar W} Z) \\
~&=& \n d \phi (Z,\bar W) -  \n d \phi (\bar W, Z)
\eean
and hence $(\n d\phi)^{1,1}$ is symmetric. This shows part (i).\\
Let us consider two Nijenhuis type complex connections $D^1$ and $D^2.$ For $k=1,2$ we write
$$\nabla^k d\phi (Z,\bar W) = \n^E_Z \Phi(\bar W) -  \Phi(D^k_Z \bar W ).$$
Hence it follows 
\bea \nabla^2 d\phi (Z,\bar W) -\nabla^1 d\phi (Z,\bar W) = \Phi(S_Z \bar W)  \label{diff_fund} \eea
for $S_Z \bar W = D^1_Z \bar W -D^2_Z \bar W.$ 
Using equation \eqref{nice_con_two} of Proposition \ref{NC_conn_prop} one has $S_Z\bar W = S_{\bar W} Z=0,$ which yields
$$ \nabla^1 d\phi (Z,\bar W) = \nabla^2 d\phi (Z,\bar W) $$
and finishes the proof of part (ii).
\epf 
\noindent
This Proposition shows, that the next Definition does not depend on the choice of the Nijenhuis type complex connection $D.$

\bd \label{pl_harm_equ_def}
Let $(M,J)$ be an almost complex manifold and $(N,h)$ a pseudo-Riemannian manifold,  then a map $ \phi \colon (M,J) \ra (N,h)$ is called {\cmssl generalized pluriharmonic} if 
it satisfies the equation \be (\nabla d\phi)^{1,1}=0, \label{pl_h_equ} \ee
where $\nabla$ is the connection on $T^*M \ot \phi^*TN$ induced by a Nijenhuis type complex connection $D$ on $(M,J)$ and the Levi-Civita connection $\n^h$ of $h.$
\ed
\noindent
One observes, that until here we did not use a (pseudo-)Riemannian metric on $(M,J).$ For an extensive study of the space 
of  Hermitian connections we refer to \cite{Gau}. \\
It is a well-known fact, that a pluriharmonic map from a (pseudo-)K\"ahler manifold into a (pseudo-)Riemannian manifold is harmonic. \\ 
Recall, that an almost complex manifold endowed with a (pseudo-)Riemannian metric is called {\cmssl almost  (pseudo-)Hermitian} if it holds $J^*g=g.$ The above mentioned result generalizes to the following proposition. Examples of this type are given in Section \ref{nksection}.
\begin{Prop}[see for instance Proposition 7 of \cite{S1}]  \label{pl_harm}
Let $(M,J,g) $ be an almost Hermitian manifold endowed with a Nijenhuis type complex connection $D$ and $(N,h)$ a pseudo-Riemannian manifold.  Moreover, suppose that the tensor $\n^g-D$ is trace-free, i.e. $\mbox{Tr}_g(\n^g_\cdot \cdot-D_\cdot \cdot )=0.$ Then a generalized pluriharmonic map  $\phi \colon (M,J) \ra (N,h)$  is harmonic.
\ep
\noindent
We  observe, that the condition $\mbox{Tr}_g(\n^g_\cdot \cdot-D_\cdot \cdot )=0$ does not depend on the choice of the Nijenhuis type complex connection $D.$
\noindent

\subsection{Associated families and symmetric targets}
In the following we fix an almost complex manifold $(M,J)$ and we put
$$\mathcal{R}_\theta= \exp(\theta  J) \in \G(\mbox{End}(TM))$$ with $\theta \in
\bR.$  An {\cmssl associated family} for $\phi\,:\, M\ra N$ is a family of maps $$\phi_\theta\,:\, M \ra N,\; \theta \in \bR,$$ such that it holds
$$ \Psi_\theta \circ d\phi_\theta = d\phi \circ \mathcal{R}_\theta,$$
for some bundle isomorphism $\Psi_\theta \colon \phi^*_\theta T N \ra \phi^* T N, \theta \in \bR,$ which is parallel
with respect to $\n^h,$ i.e. satisfies the equation $$ \Psi_\theta \circ ( \phi^*_\theta\n^h)= (\phi^*\n^h) \circ \Psi_\theta.$$
For a pseudo-Riemannian symmetric space $(N=G/K,h)$ as target manifold the next result generalizes Theorem 1 of \cite{ET}. 
\begin{Th}[Theorem 5.3 of \cite{CS1}] \label{int_family}
Let $(M,J)$ be an almost complex manifold and $(N=G/K,h)$ be a pseudo-Riemannian symmetric space. A
smooth map $\phi\colon M\ra G/K$ admits an associated family $\phi_\theta$ if and
only if it is a generalized pluriharmonic map $\phi \colon (M,J) \ra G/K$ 
and satisfies 
\be R^h(d\phi (T^{1,0}M),d\phi (T^{1,0}M))=0, \label{curv_cond_equ} \ee 
and \be \label{torsion_cond}\bar{\partial}\phi(T(V,W)) =d\phi(T(V,W)^{0,1})=0 
\ee for all $V,W \in T^{1,0}M,$
where $R^h$ is the curvature tensor of the Levi-Civita connection $\n^h$ of the metric $h$ and where $\bar \partial  \phi$ is the $(0,1)$-part of the complex linear extension $d\phi \colon TM\ot \bC \ra TN \ot \bC$ of the differential of $\phi.$
\et 
\noindent
This result motivates another generalization of the notion of a pluriharmonic map. 
\bd  \label{S1_pl_def}
A map $ \phi\colon M\ra N$ is called {\cmssl $\mathbb{S}^1$-pluriharmonic} if it admits an associated family.
\ed
\noindent
By Theorem \ref{int_family} an $\mathbb{S}^1$-pluriharmonic map is pluriharmonic in the sense of Definition \ref{pl_harm_equ_def}.
Since in the definition of an associated family we do not use the choice of a connection on $M$ this notion is independent of our choice of a Nijenhuis type connection. Moreover, as seen in Proposition \ref{nice_conn_prop} the pluriharmonic equation is independent of the choice of a Nijenhuis type complex connection. Hence we chose a connection with torsion proportional to the Nijenhuis tensor and the condition \eqref{torsion_cond} is seen to be completely determined by the almost complex structure, i.e. 
\be d\phi(N_J(V,W))=\bar{\partial}\phi(N_J(V,W)) =0, \mbox{ for all } V,W \in T^{1,0}M. \label{Nijenhuis_cond}
\ee  This condition \eqref{Nijenhuis_cond}
obviously holds true if the almost complex structure is integrable. Note, that (more directly) for an arbitrary Nijenhuis type complex connection $D$ one computes
$$ T(V,W)^{0,1} = (D_VW -D_WV -[V,W])^{0,1}=-{[V,W]}^{0,1}, V,W \in T^{1,0}M,$$
which is a multiple of the Nijenhuis tensor and yields the same observation as \eqref{Nijenhuis_cond}.
\noindent
Nonetheless as we shall see in Remark \ref{Rem_NK_full_Nij} the equation \eqref{Nijenhuis_cond} can be very restrictive. 

\subsection{Nearly K\"ahler manifolds} \label{nksection}
\noindent
A class of almost Hermitian manifolds of particular interest here is {\it nearly  K\"ahler} manifolds, 
since on the one hand one has the explicit examples of pluriharmonic maps admitting associated families given in Remark \ref{Ex_of_S1} below  and on the other hand we give examples of pluriharmonic maps without associated families in Section \ref{no_family_on_NK}. Independently, nearly  K\"ahler manifolds are of interest as one of the classes in the Gray-Hervella classification \cite{GH} and as 
an active field of research in Hermitian geometry \cites{BM,N1,N2,GDMC} and mathematical physics as backgrounds in string compactifications \cite{Strominger}.
\bd 
An almost Hermitian manifold  $(M,J,g) $ is called   {\cmssl nearly K\"ahler manifold} if it holds 
$(\n^g_X J) Y = - (\n^g_YJ)X$ for all $X,Y \in TM.$
\ed
\noindent 
On a nearly K\"ahler manifold, there exists a unique complex and metric
connection $\overline{\n}$ with totally skew-symmetric torsion \cite{FI} (Here the torsion is considered as a three-form using the metric $g$.), called the {\cmssl characteristic connection.} In fact, this connection is given by
$$ \overline{\n}_XY= \n_X^gY+ \frac{1}{2} (\n_X^gJ)JY$$
and for its torsion we have $ \overline{T}(X,Y) = (\n^g_XJ)JY.$
\bc Let $(M,J,g) $ be a nearly K\"ahler manifold. Then a generalized pluriharmonic map  $\phi:(M,J) \ra N$  is harmonic.
\ec
\pf
One easily checks (see for instance \cite{S1}), that in the nearly K\"ahler case it is 
$$ N_J(X,Y)=  4(\n_X^gJ)JY = 4 \overline{T}(X,Y)$$ 
and therefore $\overline{\n}$ is a Nijenhuis type complex connection. In addition $\overline{\n}-\n^g= \frac{1}{2}\overline{T}$ is skew-symmetric. This means, that we are in the situation of Proposition \ref{pl_harm}.
\epf
\begin{Rem} \label{Rem_NK_full_Nij} 
In real dimension six the Nijenhuis tensor, seen as the following linear map 
\bean 
\Lambda^2T^{1,0}M \ra T^{0,1}M, \\
(X,Y) \mapsto N_J(X,Y), 
\eean can be invertible.  \\
Examples of this type are given by strict nearly
K\"ahler six-manifolds, as follows for instance by combining Proposition 3.2 and Lemma
3.10 of \cite{SSH}.  \\
The most famous case is the six-sphere with its complex
structure induced by the octonions \cite{G65}. The relations
\eqref{torsion_cond} or \eqref{Nijenhuis_cond} then imply $\bar \partial \phi=0.$ 
\end{Rem}

\begin{Rem} \label{Ex_of_S1} Examples of $\mathbb{S}^1$-pluriharmonic maps from almost complex domains are given in Theorem 6 of \cite{S1}. These maps are constructed using a correspondence between $tt^*$-bundles and pluriharmonic maps from \cites{D,S}. A special class of examples for $tt^*$-bundles is provided by Levi-Civita flat nearly K\"ahler metrics, cf. Theorem 7 of the same reference. A constructive classification of flat nearly K\"ahler manifolds $(M,J,g)$ is given in \cite{CS2}. We note, that in order to obtain the non-K\"ahler or strict nearly K\"ahler examples of \cite{CS2}, one needs to admit
pseudo-Riemannian metrics. More precisely, these results give non-trivial
$\mathbb{S}^1$-pluriharmonic maps from  flat nearly pseudo-K\"ahler
manifolds $(M,J,g)$ of symmetric signature $(m,m)$ into the pseudo-Hermitian
symmetric space $SO(2m,2m)/U(m,m),$ see Section 4 of \cite{CS1}.
\end{Rem}
\noindent 
In Section \ref{no_family_on_NK} we give examples of pluriharmonic maps from nearly K\"ahler manifolds into Riemannian symmetric spaces without associated families.

\section{Non-integrable source domain and Riemannian symmetric targets}
\noindent 
The aim of this section is to give examples of pluriharmonic maps from an almost complex domain $(M,J)$ with non-integrable almost complex structure into Riemannian symmetric spaces, like for example $S^{2n}$, which admit {\bf no} associated family. In particular, we
give examples, where the source manifolds are nearly K\"ahler.  
\subsection{The twistor construction} \label{twist_con_sect}~\\
Let us first shortly recall the twistor construction  \cites{Bryant,BR,ES,EW}. One considers a $2n$-dimensional Riemannian manifold $(N,h).$ Denote by $\mathcal{J}(N)$ the bundle of almost Hermitian structures on $N,$ i.e.  $$\mathcal{J}_x(N)= \{J \in \mbox{End}(T_xN)  \, | \, J^2= -\id, \;
h(J\cdot,\cdot) =- h(\cdot,J \cdot) \}$$  and by  $\pi \, :\, \mathcal{J}(N) \ra
N, \; j_x \mapsto x$  the canonical projection. \\
The fibers of this bundle are identified with the Hermitian
symmetric space $O(2n)/U(n).$ Using the Levi-Civita connection on $(N,h)$ one
splits $$T( \mathcal{J}(N))= \mathcal H \oplus \mathcal V,$$ where 
$ \mathcal V= \ker \pi_*$  is the vertical distribution and $\mathcal H \cong
\pi^*TN$ the horizontal one. The distribution $\mathcal V$ inherits a complex structure $J^{\mathcal V}$ from $O(2n)/U(n)$
and  $\mathcal H$ acquires a complex structure $J^{\mathcal H}$ given by
$J^{\mathcal H}_j=j.$ These two can be combined to define two different almost complex structures on 
$\mathcal{J}(N):$ 
\be J_1 = J^{\mathcal V} + J^{\mathcal H} \mbox{ and } J_2 =
(-J^{\mathcal V}) + J^{\mathcal H}. \label{cx_str_def} \ee
Now we consider a second almost complex manifold $(M,J),$ then a map $\psi\, :\, M \ra \mathcal{J}(N)$ with $\pi \circ
\psi = \varphi \, :\, M \ra N$ corresponds to the choice of an almost complex structure on $\varphi^*TN,$ which is the same as the choice of  a maximally isotropic subbundle $\psi^+$ of $(1,0)$-vectors in $\varphi^*TN\ot \bC$ w.r.t. the complex bilinear extension of $h$ to $TN\ot \bC.$ The subbundle of $(0,1)$-vectors is called $\psi^-.$
The space $\mathcal{J}(N)$ endowed with these almost complex structures $J_1,J_2$ is
first introduced and studied in \cite{ES}.
\begin{Prop}[Proposition 2.3 of \cite{BR} or Theorem 5.7 of \cite{R}] \label{J2_hol_prop} ~\\
Let  $(M,J)$ be an almost complex manifold, $\psi\,: \, M \ra  \mathcal{J}(N)$
a map and set $\varphi := \pi \circ \psi.$ Then $\psi$ is holomorphic with
respect to $J_2$ if and only if 
\bea
&~&\varphi^*\n^h_{\bar Z}\;s \in \Gamma(\psi^+) \mbox{ for all }  Z \in \Gamma(T^{1,0}M)
\mbox{ and } s \in  \Gamma(\psi^+), \\
&~&d\varphi(T^{1,0}M) \subset \psi^+ 
\eea
or equivalently 
\bea
&~&\varphi^*\n^h_{ Z}\;s \in \Gamma(\psi^-) \mbox{ for all }  Z \in \Gamma(T^{1,0}M)
\mbox{ and } s \in  \Gamma(\psi^-), \\ 
&~&d\varphi(T^{0,1}M) \subset \psi^-. 
\eea
\ep
\noindent
In many geometric situations and in the applications below one has to
consider the following setting \cites{OR,BGR}: Let $\mathcal{Z}$ be a complex manifold and $p\colon \mathcal Z \ra N$  be a submersion with complex fibers and $\mathcal H^{\mathcal Z}$ a horizontal distribution for $p.$ Moreover, let $i\colon  \mathcal Z \ra \mathcal{J}(N)$ be a fiber preserving map, which is holomorphic between the fibers and maps the horizontal distribution $\mathcal{H}^{\mathcal Z}$ to $\mathcal{H}.$  We refer to  such a fibration  $p\colon \mathcal Z \ra N$ endowed with the described data $(\mathcal H^{\mathcal Z},i)$ as a twistor space. In the same manner as in Equation \eqref{cx_str_def} one obtains almost complex structures $J_1^{\mathcal Z}$  and $J_2^{\mathcal Z}$ on $\mathcal Z,$ such that the map $i$ is a holomorphic map between $(\mathcal Z,J_k^{\mathcal Z})$ and $(\mathcal{J}(N) ,J_k)$ for $k=1,2.$ If we define $\psi^+$ and $\psi^-$ in the analogous way as above we obtain the following Corollary.
\bc \label{J2_hol_cor} 
Let  $(M,J)$ be an almost complex manifold and $\mathcal Z$ a twistor space
over $(N,h),$ $\psi\,: \, M \ra  \mathcal Z$ a map and set $\varphi := p \circ \psi.$ Then $\psi$ is holomorphic with
respect to $J_2^{\mathcal Z}$ if and only if 
\bea
&~&\varphi^*\n^h_{\bar Z}\;s \in \Gamma(\psi^+) \mbox{ for all }  Z \in \Gamma(T^{1,0}M)
\mbox{ and } s \in  \Gamma(\psi^+), \label{hol_popI_equ}\\
&~&d\varphi(T^{1,0}M) \subset \psi^+ \label{hol_popII_equ}
\eea
or equivalently 
\bea
&~&\varphi^*\n^h_{ Z}\;s \in \Gamma(\psi^-) \mbox{ for all }  Z \in \Gamma(T^{1,0}M)
\mbox{ and } s \in  \Gamma(\psi^-), \label{hol_popI_equ_bis}\\
&~&d\varphi(T^{0,1}M) \subset \psi^- \label{hol_popII_equ_bis}.
\eea
\ec
\subsection{Pluriharmonic maps into symmetric targets without associated family} ~\\
In this subsection we give examples of pluriharmonic maps from almost complex manifolds into symmetric spaces $N$, 
which do not admit an associated family, since the condition \eqref{torsion_cond} of Theorem \ref{int_family} is not satisfied. In particular, if $N$ is of compact 
(as for example $S^{2n}$ or $\bC P^n$) or non-compact type, 
then the curvature condition \eqref{curv_cond_equ} of the same Theorem holds true by a Lemma of \cite{OU}, see also  p. 298 of \cite{ET}.
\bt \label{const_result}
Let $(M,J)$ be an almost complex manifold,  $(N^{2n},h)$ be a Riemannian manifold and $p \colon \mathcal Z \ra N$ be a twistor space over $(N,h),$ 
then a $J_2^{\mathcal Z}$-holomorphic map $\psi \colon M \ra \mathcal Z$ induces a pluriharmonic map $\varphi= p \circ \psi$ from $(M,J)$ to $(N,h).$
\et
\pf
Let $\psi \colon M \ra \mathcal Z$ be a $J_2^{\mathcal Z}$-holomorphic map and set $\varphi=p \circ \psi.$  For  $Z,W \in \Gamma(T^{1,0}M) $ we consider
\bea (\n_Z d\varphi) \bar W &=& \varphi^*\n^h_{ Z} (d\varphi \bar W) -d\varphi (D_Z \bar W) \label{fl_equ_11_thm} \mbox{ and }\\
(\n_{\bar W} d\varphi) Z&=& \varphi^*\n^h_{  \bar W}(d\varphi  Z) -d\varphi (D_{ \bar W}Z), \label{sl_equ_11_thm}
\eea
where $D$ is a Nijenhuis type complex connection. As $D$ is a complex connection and hence preserves the decomposition $TM\ot \bC =T^{1,0}M \oplus T^{0,1}M,$ we get using Corollary  \ref{J2_hol_cor}, that  the expression in equation \eqref{fl_equ_11_thm} lies in $\psi^-$ and the term  \eqref{sl_equ_11_thm} lies in $\psi^+.$   The fact that $(\n d\varphi)^{1,1}$ is symmetric (see Proposition \ref{nice_conn_prop}), shows  $(\n d\varphi)^{1,1}=0.$
\epf

\bc \label{cor_bspe}
Let $(N^{2n},h)$ be a Riemannian manifold and $\mathcal Z$ be a twistor space over $N,$ then $p\,:\,(\mathcal Z,J_2^{\mathcal Z}) \ra N$ is pluriharmonic.
\ec
\pf
This follows from Theorem \ref{const_result} by  setting  $\psi:=id\,:\, \mathcal Z \ra \mathcal Z.$
\epf

\begin{Ex} \label{Ex_CPn} For the following examples we can apply Theorem \ref{const_result} and show in Corollary \ref{Cor_tw_s2n_CPn}, that the condition \eqref{torsion_cond} of Theorem \ref{int_family} does not hold true.\\
 \noindent
(i) The case of the $2n$-dimensional sphere was first studied in \cite{Cal}.
In this case the twistor space is given by the homogeneous fibration
$$p\colon \mathcal Z= \mathcal{Z}(S^{2n}) = SO(2n+1)/U(n) \ra  SO(2n+1)/SO(2n) = S^{2n}$$ 
and the generic fiber is $SO(2n)/U(n),$ which is, as mentioned above, a Hermitian symmetric space and can be seen as the Grassmannian of maximal isotropic subspaces in
$\bC^{2n} \cong \bR^{2n}\otimes \bC$ (w.r.t. the bilinear extension of the Euclidean scalar product of $\bR^{2n}$). The map $i$ is given by the inclusion of $SO(2n)/U(n)$ into $O(2n)/U(n),$  which corresponds to restricting to almost complex structures compatible with a fixed orientation on $S^{2n}.$ We consider the symmetric decomposition 
$$ \mathfrak{so}(2n+1) = \mathfrak p \oplus \mathfrak{so}(2n)
\mbox{ with } \mathfrak p := \left\{ \left. \left(\begin{array}{cc}0 & v^t  \\ -v & 0_{2n \times 2n}  \end{array}
 \right) \right| v \in \bR^n\right\}  $$ 
and where $\mathfrak{so}(2n)$ is embedded as 
$$A\in \mathfrak{so}(2n) \mapsto \left(\begin{array}{cc}0 & 0  \\ 0 & A  \end{array}
 \right)\in \mathfrak{so}(2n+1) .$$ 
Let us further decompose
$$ \mathfrak{so}(2n) = \mathfrak{u}(n) \oplus \mathfrak h,$$
where $\mathfrak h$ are the Hermitian matrices. 
Setting $\mathfrak m:= \mathfrak p \oplus \mathfrak{h}$ we have the
following reductive decomposition for $\mathcal Z$ given by
  $$ \mathfrak{so}(2n+1)= \mathfrak m  \oplus\mathfrak{u}(n).$$
Let $o=e\, U(n)$ be  a canonical base point. Then $T_o \mathcal Z$ is identified with
$\mathfrak m$ and $\mathcal V_o$ corresponds to $\mathfrak h$ and $\mathcal
H_o$  to $\mathfrak p.$ In particular, from $[\mathfrak p,\mathfrak h] \subset
\mathfrak p$ it follows $[\mathcal H, \mathcal V]\subset \mathcal H.$ \\
\noindent
(ii) For $(\mathbb{C}P^n,J)$ one considers following \cites{EW,DZ} the homogeneous fibrations
$$p\colon \mathcal Z_r =\mathcal Z_r(\mathbb{C}P^n) = \frac{U(n+1)}{U(r)
  \times U(1) \times U(n-r)}  \ra  \frac{U(n+1)}{U(1) \times U(n)}
=\mathbb{C}P^n, $$ 
where $1\le r \le n$ is fixed and the symmetric decomposition
$$ \mathfrak{u}(n+1) = \mathfrak p \oplus \mathfrak{u}(1) \oplus \mathfrak{u}(n)
\mbox{ with } \mathfrak p := \left\{\left. \left(\begin{array}{cc}0 & v^h  \\ -v & 0_{n\times n}  \end{array}
 \right) \right| v \in \bC^n\right\}. $$ 
Here $  \mathfrak{u}(1) \oplus \mathfrak{u}(n)$ is embedded as 
$$(z,A)\in \mathfrak{u}(1) \oplus \mathfrak{u}(n) \mapsto \left(\begin{array}{cc}z & 0  \\ 0 & A  \end{array}
 \right)\in \mathfrak{u}(n+1) .$$ 
Moreover, let us decompose
$$ \mathfrak{u}(1) \oplus \mathfrak{u}(n)  \cong \mathfrak{u}(r) \oplus\mathfrak{u}(1) \oplus \mathfrak{u}(n-r) \oplus \mathfrak h,$$
where the $\mathfrak{u}(1)$ factors are mapped to each other and where we set $s=n-r$ and put $$\mathfrak h:=\left\{\left. \left(\begin{array}{ccc}0_{r\times r}& z_1 & A
  \\ -z_1^h&0& -z_2^h
  \\ -A^h & z_2
  & 0_{s \times s} \end{array} \right) \right| \,\, A \in M(r,s,\bC),\, z_1 \in
\bC^r,\; z_2 \in \bC^s \right\}.$$  Setting $\mathfrak m:= \mathfrak p \oplus \mathfrak{h}$ we have the following reductive decomposition for $\mathcal Z$ given by
  $$ \mathfrak{u}(n+1)= \mathfrak m  \oplus \mathfrak{u}(r) \oplus\mathfrak{u}(1) \oplus \mathfrak{u}(n-r).$$
Let $o=e\,(U(1) \times U(n)) $ be  a canonical base point. Then $T_o \mathcal Z$ is identified with $\mathfrak m$ and $\mathcal V_o$ corresponds to $\mathfrak h$ and $\mathcal
H_o$  to $\mathfrak p.$ In particular, from $[\mathfrak p,\mathfrak h] \subset
\mathfrak p$ it follows $[\mathcal H, \mathcal V]\subset \mathcal H.$ Moreover, the decomposition of $\mathfrak{p}\cong \bC^n \cong \bC^r \oplus \bC^s$ corresponds to the choice of a complex r-plane $W$ in $T_o \bC P^n$ and one defines a complex structure on $T_o \bC P^n$ by  $j= J$ on $W$ and $j=-J$ on $W^\perp \cong \bC^s.$ This yields the map $i,$ which is studied in Section 6 of \cite{OR} after identifying $\mathcal{Z}_r(\bC P^n)$ with the Grassmann bundle over $\bC P^n.$
\end{Ex}
\bc \label{Cor_tw_s2n_CPn}
There exist non-trivial (generalized) pluriharmonic maps from the twistor spaces $\left(\mathcal{Z}(S^{2n}), J_2^{\mathcal{Z}(S^{2n})}\right)$ to $S^{2n}$  
and $\left(\mathcal{Z}_r(\bC P^n),J_2^{\mathcal{Z}_r(\bC P^n)}\right),$ $1\le r \le n,$ to $\bC P^n,$ which do not admit an associated family.
\ec
\noindent 
In the proof of this Corollary we use the following Lemma relating the Nijenhuis tensors of the two complex structures on a twistor space.
\bl \label{relate_nij}
Let $\mathcal Z$ be a twistor space and denote by $N_1$ and $N_2$ the Nijenhuis tensors of the almost complex structures $J_1^{\mathcal Z}$ and $J_2^{\mathcal Z}.$ For a basic horizontal vector field $X$
and a vertical vector field $V$ one has the relation
$$ N_2 ( X, V) = N_1(X,V) +2 ([J_1V,J_1X] + J_1(\pi^{\mathcal V}[V,J_1X])).$$ In particular, the horizontal part of the difference is 
$$ \pi^{\mathcal H}\left(N_2(X,V) - N_1(X,V)\right)= 2 \, \pi^{\mathcal H}([J_1V,J_1X]),$$
where $\pi^{\mathcal H}$ and $\pi^{\mathcal V}$ are the projections to $\mathcal
H$ and $\mathcal V,$ respectively.
\el
\pf Using that for a basic vector field $X$ and a vertical vector
field $W$ the bracket $[X,W]$ is vertical\footnote{cf. for instance 9.22 of 
  \cite{Besse}} we obtain: 
\bean
N_2(X,V)&=& [J_2V,J_2X] -[V,X] -J_2[V,J_2X] -J_2[J_2V,X] \\
&=& -[J_1V,J_1X] -[V,X] -J_2[V,J_1X] +J_2[J_1V,X]\\
&=& -[J_1V,J_1X] -[V,X] -J_2[V,J_1X] -J_1[J_1V,X]\\
&=&N_1(X,V) +2 [J_1V,J_1X] -J_2[V,J_1X] + J_1[V,J_1X]\\ 
&=&N_1(X,V) + 2\left([J_1V,J_1X] + J_1 \pi^{\mathcal V}([V,J_1X]) \right) .
\eean
\epf
\pf(of Corollary \ref{Cor_tw_s2n_CPn})
As in Corollary \ref{cor_bspe} one can take the identity for the map $\psi.$
 In both cases it is known that $J_1$ is integrable.\footnote{ For $S^{2n}$ the integrability of $J_1$ follows, since $S^{2n}$ is conformally flat. Moreover, for an almost Hermitian manifold $N$ the integrability of $J_1$ on the Grassmann bundle of $k$-planes is equivalent to the vanishing of the Bochner curvature, cf. Section 6 of \cite{OR}. Since $\bC P^n$ has constant holomorphic sectional curvature, it is Bochner flat and $J_1$ is integrable.}
For the Nijenhuis tensor $N_2$ of $J_2$ it follows from Lemma \ref{relate_nij} and the vanishing of $N_1,$ that one has $dp\circ N_2\ne 0,$ where we used the relation $[\mathcal H,\mathcal V] \subset \mathcal H$ given in Example \ref{Ex_CPn}. Therefore the necessary condition \eqref{torsion_cond} is not satisfied.
\epf

\subsection{Pluriharmonic maps from nearly K\"ahler manifolds without associated family} \label{no_family_on_NK} ~\\
In this section we are interested in nearly K\"ahler structures $(M,J,g)$ on twistor
spaces.  Further we suppose, that $M$ is a  strict nearly K\"ahler manifold and has
special algebraic torsion in the sense of the next Definition \cite{N1}.
\bd
A strict nearly K\"ahler manifold $(M,J,g)$ has {\cmssl special algebraic torsion}, if there exists an orthogonal decomposition $TM=\mathcal H \oplus \mathcal V,$ which is parallel for the characteristic connection $\overline{\n},$ such that $\mathcal H$ and $\mathcal V$ are invariant under $J$ and it holds $\overline{T}(\mathcal V,\mathcal V)=0$ and
$\overline{T}(\mathcal H,\mathcal H)= \mathcal V.$ 
\ed
\noindent
Additionally assume, that $M$ is connected, complete and irreducible (as Riemannian manifold), then there exist a Riemannian manifold $(N,g_N)$ and a
Riemannian submersion with totally geodesic fibers $p\colon M \ra N,$ called
the {\cmssl canonical} fibration, such that the fibers with the induced
Riemannian metric and almost complex structure are Hermitian symmetric
spaces and such that the vertical bundle of $p$ coincides with $\mathcal V.$ For more details we refer to Section 4 of \cite{N1} or Section 1 of \cite{GDMC}. \\
Denote by $g_{\mathcal V},$ $g_{\mathcal H},$ $J_{\mathcal V}$ and
$J_{\mathcal H}$ the restrictions of the nearly K\"ahler data $J$ and $g$  to $\mathcal V$ and $\mathcal H.$ Then $\hat g:= 2g_{\mathcal V} \oplus g_{\mathcal H}$ and $\hat J$ with $\hat J_{\mathcal V} = -J_{\mathcal V}$ and $\hat J_{\mathcal H} = J_{\mathcal H}$
give a K\"ahler structure $(\hat J,\hat g)$ on $M.$ By the classification of
\cite{N1} (see also \cite{GDMC}) $p$ is a twistor fibration and one has the following information:
\begin{itemize}
\item[(i)] $M$ is homogeneous if and only if $N$ is a symmetric space,
\item[(ii)] $M$ is non-homogeneous if and only if $N$ is a non-symmetric
  positive quaternionic K\"ahler manifold $N^{4k}$ and $M$ is its twistor space.
\end{itemize}
Let us note, that $J$ and $\hat J$ correspond to $J_1$ and $J_2,$ respectively.
The map $i$ is in both cases (i) and (ii) the inclusion of the fibers of $p\colon M \ra N$ of the projection $p$ into $M,$ see Proposition 4.2 of \cite{N1}. \\
For completeness sake we shortly recall, that a quaternionic K\"ahler manifold is a Riemannian manifold $(N^{4k},g)$ of dimension $4k$ with holonomy $Sp(k)\cdot Sp(1).$ 
This means, that there exists a rank 3 sub-bundle $Q$ of $\mbox{End}(TN),$ which is preserved by the Levi-Civita connection of $g$ and is locally generated by three anti-commuting endomorphisms-fields $j_1,j_2,j_3,$ that satisfy $j^2_1=j^2_2=-\id$ and $j_3= j_1 j_2 =-j_2 j_1.$ The condition that $Q$ is preserved by the Levi-Civita connection $\n^g$ is
in a given standard local basis $\{j_\alpha\}_{\alpha=1}^3$ of $Q$  equivalent to
the equations 
\be
 \n^g_X j_\alpha = -\theta_\gamma (X) j_\beta + \theta_\beta(X) j_\gamma, \mbox{ for } X \in TN, 
\label{Q_parallel}\ee
where $\alpha,\beta,\gamma$ is a cyclic permutation of $1,2,3$ and $\{\theta_\alpha \}_{\alpha=1}^3$
are local one-forms. One defines the twistor space as the following $S^2$-sphere bundle 
$$ \mathcal Z^{QK}:= \{ A \in Q \, | \, A^2=-\id \}.$$
The fibers carry the structure of a Hermitian symmetric space and hence that of a K\"ahler manifold and one applies a similar construction as given in section \ref{twist_con_sect} to define two almost complex structures $J_1,J_2$ on $ \mathcal Z^{QK},$ such that there is a holomorphic inclusion into $\mathcal{J}(N)$ which is induced by the inclusion of the corresponding endomorphisms. Moreover,
one obtains an almost Hermitian metric by 
$$g_z = \frac{1}{\nu} g^{S^2} \oplus g_{\pi (z)}, \mbox{ for } z \in \mathcal Z^{QK},$$ where $\nu=\frac{\mbox{scal}}{4k(k+2)}$ is the reduced scalar curvature.
This metric is a K\"ahler-Einstein metric by Theorem 6.1 of \cite{Sal} or Theorem 14.80 of \cite{Besse}. 
The corresponding nearly K\"ahler examples $(\hat J, \hat g)$ of this type can also be characterized by a holonomic condition \cites{BM,N2}. Combining the above classification information with Theorem \ref{const_result} we
obtain the following examples.
\bp \label{bsp_prop_tw_space}
Let $(\mathcal Z^{QK},g,J)$ be the twistor space of a positive
quaternionic K\"ahler manifold endowed with its canonical nearly K\"ahler structure. \\ Then the projection $p \colon \mathcal Z^{QK} \ra N^{4k}$ is a (generalized)
pluriharmonic map, which does not admit an associated family.
\ep 
\pf
As above we use Theorem \ref{const_result} with the identity map as $\psi \colon \mathcal Z^{QK} \ra \mathcal Z^{QK},$ then the projection is pluriharmonic.\\
 Let us recall, that we are in the situation of special algebraic torsion, i.e. we have a $\overline{\n}$-parallel orthogonal and $J$-invariant decomposition $TM=\mathcal H \oplus \mathcal V.$ \\
{\bf Claim:} It holds  $\{0\} \ne \overline{T}(\mathcal H, \mathcal V) \subset
\mathcal H.$  \\ Firstly, we have $0=g(\overline{T}(U,V),X) =g(\overline{T}(X,U),V), \mbox{ for } U,V \in \mathcal V \mbox{ and } X \in \mathcal H,$ i.e. $\overline{T}(X,U)$ lies in $\mathcal H$ for $U \in \mathcal V \mbox{ and } X \in \mathcal H.$ \\ Secondly, we show that $$W:=\mbox{span} \{ \overline{T}(X,U) |  U \in \mathcal V, X \in \mathcal H \} \ne \{0\}.$$  Let us suppose $W=\{0\},$ then it follows
$0=g(\overline{T}(X,U),Y)=-g(\overline{T}(X,Y),U), \forall U \in \mathcal V, X,Y \in \mathcal H, $
which contradicts the special algebraic torsion condition
$\overline{T}(\mathcal H,\mathcal H)= \mathcal V$  and finishes the  proof of the claim.\\
Since the Nijenhuis tensor
is a multiple of $\overline{T},$ we conclude that the necessary condition
\eqref{torsion_cond} of Theorem \ref{int_family} is not valid.
\epf 
\noindent
In Section 8 of \cite{GDMC} detailed lists of homogeneous nearly K\"ahler
manifolds with special algebraic torsion and the according canonical
fibrations combining the results of \cites{G72,N1} are given. These give
further examples of  (generalized) pluriharmonic maps without associated families.
\bp
Let $(M,J,g)$ be a connected, complete and irreducible (as Riemannian manifold)
 homogeneous strict nearly K\"ahler manifold with special algebraic
torsion, then the canonical fibration $p\colon M \ra N$ is  a (generalized)
pluriharmonic map to the symmetric space $N,$ which does not admit an associated family.
\ep
\pf The proof goes along the same line as that of Proposition
\ref{bsp_prop_tw_space}. \epf 
\noindent
Let us recall, that $S^4$ and $\mathbb{C}P^2$ are the only four-dimensional compact positive quaternionic K\"ahler manifolds 
and their twistor spaces, namely $\mathbb{C}P^3$ and the flag manifold $F(1,2),$ are by results of \cite{Bu} the only (simply connected) homogeneous examples with special algebraic torsion in dimension six. 
\bc
There exist non-trivial (generalized) pluriharmonic maps from the
six-dimensional nearly K\"ahler manifolds $\mathcal{J}(S^4)$ and
$\mathcal{J}(\mathbb{C}P^2),$ resp. endowed with the above given nearly K\"ahler structure to the Riemannian symmetric spaces $S^4$ and $\bC P^2,$ resp., which do not admit an associated family.
\ec 
\section{Integrable source and pseudo-Riemannian targets}
\noindent 
In this section we consider a smooth map $f \colon M \ra N=G/H$ from a complex
manifold $(M,J)$ into a reductive homogeneous space $N=G/H$  with canonical projection $\pi \colon G \ra G/H.$ Further let 
\bea \mathfrak g = \mathfrak h \oplus \mathfrak m \label{hom_dec} \eea
be the associated  reductive decomposition. Denote by $\theta \colon TG \ra \mathfrak g$ the left Maurer-Cartan form of $G,$ which splits according to Equation \eqref{hom_dec} into its components 
$$\theta= \theta_{\mathfrak{m}} \oplus \theta_{\mathfrak{h}}.$$
A map  $F \colon M \ra G$  is called framing\footnote{For a contractible manifold such a framing always exists. This fact is not needed in the present text, as we shall always start with some given framing.} of $f$ provided that it holds 
$f:= \pi \circ F.$ For a framing $F$ let $\alpha:= F^*\theta$ and decompose $\alpha$ with respect to the direct sum \eqref{hom_dec} into  
$$\alpha= \alpha_{\mathfrak{m}} \oplus \alpha_{\mathfrak{h}}.$$
The tangent space at $o = eH$ of $G/H$ is identified with $\mathfrak
g/\mathfrak h$ and as $G	/H$ is reductive with $\mathfrak m$ via
$$ {\mathfrak g}/{\mathfrak h} \cong \mathfrak{m} \ni X \mapsto \left. \frac{d}{dt}\right|_{t=0} exp(tX)\cdot o.$$
This provides an identification of $TN$ with $G \times_H \mathfrak m.$
Moreover, there exists a natural inclusion of $G \times_H \mathfrak m$ into $G \times_H
\mathfrak g,$ which itself can be canonically identified with $N \times
\mathfrak g.$ Combining the above information we get an identification 
$$ \beta \colon TN \ra N \times \mathfrak g$$
of $TN$ with a subbundle of $N \times \mathfrak g$ and we may interpret $\beta$ as a $\mathfrak
g$-valued 1-form on $N.$\\
 Let us now additionally suppose, that $(G/H,h)$ is naturally reductive, i.e. $h$ is a $G$-invariant pseudo-Riemannian metric on $G/H,$  such that one has 
$$ \langle [X,Y]_{\mathfrak m}, Z \rangle = \langle X, [Y,Z]_{\mathfrak m}\rangle, \mbox{   for } X,Y,Z \in \mathfrak m,$$
where $\langle \cdot ,\cdot \rangle$ is the pseudo-Euclidean scalar product on $\mathfrak m$ corresponding to $h.$
\begin{Prop}[Chapter 1 of \cite{BR}, Section 1 of \cite{H98}] \label{con_lc_conn_propo}
The Levi-Civita connection $\n^h$ of the naturally reductive homogeneous space $(N=G/H,h)$ is given by 
$$ \beta(\n^h_XY)= X(\beta(Y)) - [\beta(X),\beta(Y)] +\frac{1}{2}
P[\beta(X),\beta(Y)], \; X,Y \in \Gamma(TN),$$
where $P$ denotes the projection onto the tangent bundle 
$P \colon N \times \mathfrak g \ra G \times_H \mathfrak m$
and the canonical connection is given by
$$ \beta(\n^{can}_XY)= X(\beta(Y)) - [\beta(X),\beta(Y)],\; X,Y \in \Gamma(TN).$$
\ep
\noindent
From this Proposition one obtains the observation, that
if $G/H$ is symmetric, then $\n^{can}$ and $\n^h$ coincide (see for example
the Remark on p. 10 of \cite{BR}). Moreover, we have the relations 
$$ (\pi^* \beta)_g= \mbox{Ad}(g)(\theta_{\mathfrak m})_g \mbox{ and } f^*\beta
= \mbox{Ad}(F) \alpha_{\mathfrak m}.$$
\bp \label{pl_harm_equ_in_framing}
Let $f \colon (M,J) \ra (G/H,h)$ be a smooth map from a complex manifold $(M,J)$ into a pseudo-Riemannian symmetric space $(G/H,h)$ and let $F \colon M \ra G$ be a framing of $f.$ The map $f$ is pluriharmonic if and only if 
$$ \bar \partial \alpha_{\mathfrak m}'
+\left[\alpha_{\mathfrak{h}}'' \wedge \alpha_{\mathfrak m}' \right] =0,$$
where $\alpha_{\mathfrak m}'$ is the $(1,0)$-part of $\alpha_{\mathfrak m}$ 
and  $\alpha_{\mathfrak h}''$ is the $(0,1)$-part of $\alpha_{\mathfrak
  h}.$\ep 
\pf
Using Proposition \ref{con_lc_conn_propo} and the identity \be \beta(\partial f)=(f^*\beta)' \label{help_id} \ee one has for the $(1,1)$-part of $\n df$ 
$$\beta\left({\n}'' \partial f\right)= \bar \partial \beta(\partial f)
-[(f	^*\beta)''\wedge (f	^*\beta)'],$$
where $\n$ is the connection induced by $\n^{can}=\n^h$ and by a Nijenhuis type complex connection $D.$  Applying the identity \eqref{help_id} and $(f^*\beta)'= \mbox{Ad}(F)\,\alpha'_{\mathfrak m}$ we obtain
\begin{eqnarray*}
\bar \partial \beta(\partial f)&=& \bar \partial (\mbox{Ad}(F)\, \alpha_\mathfrak
m') \\
&=&\mbox{Ad}(F)\, \left( \bar \partial \alpha_{\mathfrak m}' +
  \left[\alpha'' \wedge \alpha'_{\mathfrak m}\right] \right) \\
&=&\mbox{Ad}(F) \left( \bar \partial \alpha_{\mathfrak m}' +
  \left[\alpha_{\mathfrak h}''\wedge \alpha_{\mathfrak m}'\right]+\left[\alpha_{\mathfrak m}''\wedge \alpha_{\mathfrak m}'\right] \right). 
\end{eqnarray*}
Further it is
$$[(f^*\beta)''\wedge (f^*\beta)']= [(\mbox{Ad}(F)\, \alpha_\mathfrak
  m'')\wedge (\mbox{Ad}(F)\, \alpha_\mathfrak m')]=  \mbox{Ad}(F)\, [ \alpha_\mathfrak
  m'' \wedge  \alpha_\mathfrak m'] $$
and therefore we get
$$\beta\left({\n}'' \partial f \right)= \mbox{Ad}(F)\, \left( \bar \partial \alpha_{\mathfrak m}'
+\left[\alpha_{\mathfrak{h}}'' \wedge \alpha_{\mathfrak m}' \right]\right).$$
Since $\beta  \colon TN \ra N \times \mathfrak g$ is injective, this equation finishes the proof.
\epf 

\noindent
After these preparations the aim is to construct examples of pluriharmonic maps where the curvature condition \eqref{curv_cond_equ} of Theorem \ref{int_family} is {\bf not} satisfied. \\
For a non-compact semi-simple Lie algebra $\mathfrak h$ we consider its complexification $\mathfrak
g:=\mathfrak h^\bC$ and the symmetric decomposition  \be \mathfrak g= \mathfrak h^\bC =\mathfrak h
\oplus \sqrt{-1}\mathfrak h= \mathfrak{p}\oplus \mathfrak{h}, \label{symdec}  \mbox{ with } \mathfrak{p}= \sqrt{-1}\mathfrak h. \ee
Denote by $H$ and $G$  Lie groups with Lie algebras $\mathfrak h$ and $\mathfrak g=\mathfrak h^\bC,$ such that $H$ is a closed subgroup in $G.$  An invariant (pseudo-)Euclidean scalar product $\langle \cdot,\cdot \rangle$ on $\mathfrak h$ can be extended by $\mbox{Re}(\langle \cdot,\bar \cdot \rangle)$ to $\mathfrak g=\mathfrak h^\bC$  and endows  $N=G/H$ with the structure of a pseudo-Riemannian symmetric space. Since $\mathfrak h$ is semi-simple,  a possible choice of $\langle \cdot, \cdot
\rangle$ is given by the Killing form. \\
Composing a map $F \colon M \ra G$ with the canonical projection $\pi\colon G \ra G /H$ yields a map $f \colon M \ra G /H$ with framing $F.$  For a general map $F\colon M \ra G$ the integrability conditions\footnote{i.e. the obstruction to the local existence of a map $F\colon M \ra G$ with $\alpha=F^*\theta$}  in terms of $\alpha=F^*\theta$ are given by the Maurer-Cartan equations
$$ d\alpha + \frac{1}{2} [\alpha \wedge \alpha] =0.$$
Decomposing with respect to \eqref{symdec} one gets  
\bea
&~& d\alpha_{\mathfrak p} +  [\alpha_{\mathfrak h} \wedge \alpha_{\mathfrak
    p}] =0, \label{MC_1} \\
&~& d\alpha_{\mathfrak h} + \frac{1}{2} \left([\alpha_{\mathfrak h} \wedge
   \alpha_{\mathfrak h}] +[\alpha_{\mathfrak p} \wedge   \alpha_{\mathfrak p}]\right) =0,
\eea
where $\alpha_{\mathfrak h}$ and $\alpha_{\mathfrak p}$ are the $\mathfrak h$-
and the $\mathfrak p$-part of the $\mathfrak g$-valued one-form $\alpha.$ \\
\noindent 
In the following we consider  a Borel subgroup $(M,J)=(B,J)$ of $G$ with Lie algebra
$\mathfrak b$ and let $F=\iota \colon B  \hookrightarrow G$ be its
inclusion into $G,$ which is a holomorphic map\footnote{see for example p. 506 of \cite{Kn}}.\\
More precisely, let $\mathfrak a$ be a Cartan subalgebra of $\mathfrak g =
\mathfrak h^\bC$ and $\Delta^+$ be a positive root system, then any Borel
subalgebra is conjugate to
$$ \mathfrak b_0=\mathfrak a  \oplus \sum_{\alpha \in  \Delta^+} \mathfrak
h^\alpha = \mathfrak a \oplus \mathfrak n^+ \mbox{ with } \mathfrak n^\pm =
\sum_{\pm \alpha \in \Delta^+} \mathfrak h^\alpha$$
and one directly checks $ [\mathfrak b_0, \mathfrak b_0] = \mathfrak n^+$ and
consequently  $ [\mathfrak b, \mathfrak b] = \mathfrak n$ for some nilpotent
subalgebra $\mathfrak n.$ Since $\mathfrak h$ is semi-simple, it follows, that
its complexification is semi-simple and therefore $\mathfrak n$ is non-empty. 
\\
\noindent
The equations \eqref{MC_1} and the holomorphicity of $F\colon (B,J) \ra G$ imply
$$  d(\alpha_{\mathfrak p} \circ J ) +  [\alpha_{\mathfrak h} \wedge  (\alpha_{\mathfrak p}\circ J)] =0.$$
This equation and equation \eqref{MC_1} show using Proposition \ref{pl_harm_equ_in_framing},
that $f=\pi\circ F$ is pluriharmonic. By our assumption, that $\mathfrak h$ is a
non-compact real form of $\mathfrak g,$ $f$ is a non-trivial map.
\noindent
As $B$ is a Borel subgroup, it follows that the image of the curvature operator 
$$\mathcal R^h(dF (T^{1,0}B), dF(T^{1,0}B))$$ corresponds to $[\mathfrak b, \mathfrak
  b] = \mathfrak n,$ where $\mathfrak n$ is a nilpotent subalgebra in
$\mathfrak g.$ Therefore $\mathcal R^h$ takes values in a non-trivial isotropic subspace of
$\mathfrak g$  and hence the map $F=\iota$ violates the curvature condition \eqref{curv_cond_equ}. 
Summarizing we obtain the following result.

\bt
Let $\mathfrak h$ be a non-compact semi-simple Lie algebra and $\mathfrak g$ its
complexification. Denote by $(G/H,h)$ a (locally) pseudo-Riemannian symmetric space 
with Cartan decomposition $\mathfrak g = \mathfrak h + \sqrt{-1} \mathfrak h$
and by $\pi \colon G \ra G/H$ the canonical projection. Fix a Borel subgroup $B$
of $G$ and denote by $F\colon B \ra G$ its inclusion into $G,$ then the map $f:=
\pi \circ F$ is a pluriharmonic map, which  violates the curvature condition \eqref{curv_cond_equ} of Theorem \ref{int_family}. In particular, the map $F$ does not admit an associated family.
\et
\begin{Ex}
Let us consider $\mathfrak h = \mathfrak{su}(p,q),$ then one has $\mathfrak g= \mathfrak{sl}(p+q,\bC).$ A possible choice for $\mathfrak{n}$ is given by the space of strictly upper triangular matrices and $\mathfrak b$ is given by the intersection of the space of  upper triangular matrices with  $\mathfrak{sl}(p+q,\bC).$ \\
In order to compute the pluriharmonic map $f$ explicitely, let us determine
the Cartan embedding (see for example Ch. 3 of \cite{CE}) of $SL(p+q,\bC)/SU(p,q)$ into $SL(p+q,\bC).$ Denote by $I_{p,q}$
the matrix
$$ I_{p,q}=  \left(\begin{array}{cc} \id_p & 0  \\  0 & -\id_q \end{array} \right)$$
representing the Hermitian inner product on $\bC^{p,q}.$ The Cartan involution on
$\mathfrak g$ is given by
\begin{align*} 
\sigma \colon  \mathfrak{sl}(p+q,\bC) &\ra \mathfrak{sl}(p+q,\bC),  \\
 A \quad &\mapsto \quad I_{p,q}\, \bar{A}^t \, I_{p,q}
\end{align*}
and is the differential of the map
\begin{align*} 
\Sigma \colon SL(p+q,\bC) &\ra SL(p+q,\bC), \\
  g ~~\quad &\mapsto \quad I_{p,q}\, \bar{g}^t \, I_{p,q}.
\end{align*} 
Using this information we obtain the following explicit representation of the
Cartan embedding $\Psi$ of  $SL(p+q,\bC)/SU(p,q)$  into $SL(p+q,\bC):$ 
\begin{align*} 
\Psi \colon SL(p+q,\bC) /SU(p,q) &\ra SL(p+q,\bC), \\
  gH  \quad &\mapsto \quad g \Sigma\left(g^{-1}\right) = g \, I_{p,q}\,
\left(\bar{g}^t\right)^{-1} \, I_{p,q}= g \, I_{p,q}\,
\left(I_{p,q}\,\bar{g}^t\right)^{-1}.
\end{align*}
Let us recall, that $\Psi$ is a totally geodesic embedding. Therefore $f$ is pluriharmonic if and only if $\Psi \circ f$ is pluriharmonic. An element $b \in B$ can be written as 
\bean
 b= \left(
    \begin{array}{cc}
     A_1 &  C\\
     0 & A_2 
    \end{array}  \right),
\eean
where $A_1,A_2$ are upper triangular matrices and $C$ is a complex $p\times q$ matrix. In this representation one computes
$$
b \cdot I_{p,q}= \left(
    \begin{array}{cc}
     A_1 &  -C\\
     0 & -A_2 
    \end{array}  \right), \quad  \bar{b}^t = \left(
    \begin{array}{cc}
     \bar{A}_1^t &  0\\
     \bar{C}^t & \bar{A}_2^t 
    \end{array}  \right) \mbox{ and }  I_{p,q} \cdot \bar b^t = \left(
    \begin{array}{cc}
     \bar A_1^t &  0\\
     -\bar C^t & -\bar A_2^t 
    \end{array}  \right).
$$
Summarizing this shows the following formula for the pluriharmonic map $f$ composed with
the Cartan embedding
$$(\Psi \circ f)(b)=  \left(
    \begin{array}{cc}
     A_1 &  -C\\
     0 & -A_2 
    \end{array}  \right)  \left(
    \begin{array}{cc}
     \bar{A}_1^t &  0\\
     -\bar C^t & -\bar{A}_2^t 
    \end{array}  \right)^{-1}. $$
\noindent 
An example of this type is the {\it ``cooling tower``} (one-sheeted hyperboloid), where $\mathfrak
g=\mathfrak h^\bC= \mathfrak{sl}(2,\bC) \cong \mathfrak{so}(3,1)$ and
$\mathfrak h= \mathfrak{sl}(2,\bR) \cong \mathfrak{so}(2,1).$ In this case one
has
\bean
  B \ni b=\left(\begin{array}{cc}
     a_{1} &  c\\
     0 & a_{2} 
    \end{array}  \right), \mbox{ for } a_1,\, a_2,\, c \in \bC.
\eean
Using $\det(b)= a_{1}a_{2}=1$ we obtain the subsequent explicit form for the pluriharmonic map $\Psi \circ f$
\bean (\Psi \circ f)(b)&=&  \left(
    \begin{array}{cc}
     a_{1} &  -c\\
     0 & -a_{2} 
    \end{array}  \right)  \left(
    \begin{array}{cc}
     \bar{a}_{1} &  0\\
     -\bar{c} & -\bar{a}_{2}
    \end{array}  \right)^{-1} \\ 
&=& -\frac{1}{ \bar{a}_{1} \bar{a}_{2}} \left(
    \begin{array}{cc}
     a_{1} &  -c\\
     0 & -a_{2} 
    \end{array}  \right)  \left(
    \begin{array}{cc}
     -\bar{a}_{2} &  0 \\
     \bar{c} & \bar{a}_{1}
    \end{array}  \right) \\
&=& -\frac{1}{ \bar{a}_{1} \bar{a}_{2}} \left(\begin{array}{cc}
     -a_{1} \bar{a}_{2}-|c|^2& - \bar{a}_{1} \, c\\
     -a_{2}\, \bar{c} & -a_{2} \bar{a}_{1}
    \end{array}  \right)=\left(\begin{array}{cc}
     a_{1} \bar{a}_{2} +|c |^2& c \bar{a}_{1}\\
     a_{2} \bar{c} & a_{2} \bar{a}_{1}
    \end{array}  \right).
\eean
\end{Ex}

\end{document}